

\input amstex

\documentstyle{amsppt}
\magnification = \magstep1
\hsize = 6.25 truein
\vsize = 22 truecm
\baselineskip .22in
\define\Lap{\varDelta}
\define\del{\partial}
\define\a{\alpha}
\predefine\b{\barunder}
\redefine\b{\beta}
\define\ga{\gamma}
\predefine\d{\dotunder}
\redefine\d{\delta}
\define\th{\theta}

\define\e{\epsilon}

\predefine\o{\orsted}
\redefine\o{\omega}
\define\Sig{\Sigma}
\define\Lam{\Lambda}
\define\lam{\lambda}

\define\cin{\Cal C^{\infty}}
\define\RR{\Bbb R}
\define\SS{\Bbb S}
\define\HH{\Bbb H}
\define\TT{\Bbb T}
\define\h{\text{\bf H}}
\define\A{\text{\bf A}}

\define\ML{\Cal M_{\Lambda}}
\define\Mk{\Cal M_{k}}

\define\ind{\text{ind\,}}
\define\relind{\text{rel-ind\,}}

\define\BL{\Cal B}

\NoRunningHeads
\topmatter
\title The Moduli Space of Complete Embedded Constant Mean Curvature
Surfaces \endtitle
\author
Rob Kusner${}^{(\dag)}$, Rafe Mazzeo${}^{(\ddag)}$ and
Daniel Pollack${}^{(\star)}$
\endauthor
\abstract
We examine the space of surfaces in $\RR^{3}$ which are complete,
properly embedded and have nonzero constant mean curvature.
These surfaces are noncompact
provided we exclude the case of the round sphere.  We prove that the space
$\Mk$ of all such surfaces with $k$ ends (where surfaces are identified
if they differ by an isometry of $\RR^{3}$) is locally a real analytic
variety.  When the linearization of the quasilinear elliptic
equation specifying mean curvature equal to one has no $L^2-$nullspace
we prove that $\Mk$ is locally the quotient of a real analytic manifold of
dimension $3k-6$ by a finite group (i\.e\. a real analytic orbifold),
for $k\geq 3$.  This finite group is the isotropy subgroup of
the surface in the group of Euclidean motions.
It is of interest to note that the dimension of $\Mk$
is independent of the topology of the underlying punctured Riemann
surface to which $\Sig$ is conformally equivalent.  These results
also apply to hypersurfaces of $\HH^{n+1}$ with nonzero constant
mean curvature greater than that of a horosphere and whose ends
are cylindrically bounded.
\endabstract
\affil
Mathematical Sciences Research Institute ($\dag$, $\ddag$ and $\star$),
Stanford University ($\ddag$), University of Chicago ($\star$),
and University of Massachusetts at Amherst ($\dag$)
\endaffil
\thanks Research supported in part by
($\dag$) NSF grant \# DMS9404278 and an NSF Postdoctoral Fellowship,
($\ddag$) NSF Young Investigator Award, a Sloan Foundation Postdoctoral
Fellowship and NSF grant \# DMS9303236, and
($\star$)  NSF grant \# DMS9022140 and an NSF Postdoctoral Fellowship.
\endthanks
\endtopmatter
\document

\specialhead I. Introduction \endspecialhead

A fundamental object of study in classical differential geometry is the class
of hypersurfaces with constant mean curvature in the simply connected
space forms.  We focus
on the class of those surfaces in $\RR^{3}$ (and hypersurfaces in $\HH^{n+1}$
with an additional hypothesis) which, in addition to having constant,
nonzero mean curvature, are embedded and complete.
We will always use embedded to mean properly embedded.
Our aim is to completely understand the moduli space of all such surfaces.
In this note we describe the local structure of this space.

The main body of this paper treats the case of surfaces in $\RR^3$;
the applications to hypersurfaces in $\HH^{n+1}$ are discussed in \S 5.
We take the mean curvature of a submanifold to be the sum of the principal
curvatures rather than the average, so that the unit sphere in $\RR^{3}$
has constant mean curvature $\h\equiv 2$.  A compact, embedded surface of
constant mean curvature is necessarily a round sphere, thus the surfaces
considered here are all noncompact.  We define the ``ends'', $E_{j}$, of
an embedded surface $\Sig$ with finite topology to be the noncompact connected
components of the surface near infinity i.e.
$$
\Sig\cap(\RR^{3}\setminus B_{R}(0)) = \cup_{j=1}^{k}E_{j},
$$
where $B_{R}(0)$ denotes the ball of radius $R$ about the origin and
$R$ is chosen sufficiently large so that the number $k$ is constant
for all $R^{\prime}>R$.
For surfaces of constant nonzero mean curvature $\h$ in $\RR^{3}$, the
sign of ${\text{\bf H}}$ and its particular value may be changed by a
reversal of orientation and a homothety respectively.  Thus, we always
normalize our surfaces to have mean curvature $\h\equiv 1$.
The canonical example of such a surface is the unit cylinder.
In 1841, C.~Delaunay [1] discovered a one-parameter family of embedded
constant mean curvature (CMC) surfaces of revolution.  These surfaces
are periodic and interpolate between the unit cylinder and the
singular surface formed by a string of spheres of radius 2, each
tangent to the next along a fixed axis. In particular, he established
that every CMC surface of revolution was necessarily one of these
``Delaunay surfaces''.  Examples of embedded, complete CMC surfaces
with more complicated topology were finally found in 1987 by
N. Kapouleas [4].  Kapouleas produced his surfaces by first constructing
surfaces which are approximately CMC and then perturbing them to
nearby CMC surfaces by solving the partial differential equation which
prescribes the mean curvature to 1. This construction yields
embedded, complete surfaces with arbitrary genus and $k$ ends,
$k\geq3$.  Kapouleas also constructed, both compact and noncompact,
immersed CMC surfaces.  K. Grosse-Brauckmann [2] has used conjugate
surface constructions to obtain families of symmetric embedded
complete CMC surfaces. In particular he has shown that the genus zero
surfaces with $k$ ends possessing the maximal symmetry ($k$-fold
dihedral symmetry, interchanging the ends and reflecting across a plane)
lie in a one parameter family, realizing each asymptotic
Delaunay surface twice except for an extremal asymptotic surface
$D_{\bar{\e}}, {\bar{\e}}<1$, which is uniquely realized.
$\bar{\e}$ is determined by the requirement that the ratio of the
$\bar{\e}$ to the maximum `bulge' (see below) is equal to
$1/(k-1)$.

Results of Meeks [9], and Korevaar-Kusner-Solomon [7]  establish that a
CMC surface with at most two ends is necessarily Delaunay.
Furthermore, the main result of [7] implies that every embedded, complete
CMC surface is ``asymptotically Delaunay''.  This result is a central
ingredient in our analysis of the space of all CMC surfaces and is
recalled more carefully in the next section.

We define $\Mk$ to be the space of all complete, embedded CMC surfaces in
$\RR^{3}$ with $k$ ends, where two surfaces are considered equivalent if
they differ by a rigid motion of $\RR^{3}$.
In [5], Korevaar and Kusner examined the degree to which general
elements $\Sig\in\Mk$ `look like' the surfaces constructed by Kapouleas.
In particular, they showed that any $\Sig\in\Mk$ is contained in a regular
neighborhood of a piecewise linear graph in $\RR^{3}$.  This neighborhood
is given by a union of solid, half-infinite cylinders of radius 6, cylindrical
segments of radius 6 and solid balls of radius 21; the number of each type
is bounded by the topology of $\Sig$.  Consequently they established
a priori area and curvature bounds for $\Sig\in\Mk$.  This led to a local
compactness result and suggested that $\Mk$ is a finite dimensional
real analytic variety, which is what we establish here.

Let $\Sig\hookrightarrow\RR^{3}$ be a complete, embedded CMC surface
with $k$ ends, with unit outward normal $\nu$.
One way to obtain nearby surfaces is as follows.
Given a function $\phi\in C^{\infty}(\Sig)$, we define
a new surface $\Sig_{\phi}$ by
$$
\Sig_{\phi}=\{x+\phi(x)\nu(x): x\in\Sig\}.
$$
$\Sig_{\phi}$ will  also be embedded provided that $\phi$ is
sufficiently small in $C^2-$norm.   We need to consider
a slightly larger class of variations of $\Sigma$ than just those
arising from normal graphs in order to allow for variations of the
directions of the ends of $\Sig$. These other variations are discussed
more carefully below.
The mean curvature, $\h_{\phi}$, of $\Sig_{\phi}$ may be expressed in terms of
$\phi$ and $\Sig\hookrightarrow\RR^{3}$ as follows
$$
N(\phi)\equiv 2(\h_{\phi}-1) = \Lap{\phi}+|\A|^{2}\phi+
Q(\phi,\nabla\phi,\nabla^{2}\phi), \tag 1.1
$$
where $\Lap$ denotes the Laplace-Beltrami operator on $\Sig$,
$|\A|^{2}$ is the squared norm of the second fundamental form $\A$ of
$\Sig$, and $Q(\phi,\nabla\phi,\nabla^{2}\phi)$ is the quadratically vanishing
nonlinearity. An explicit calculation of $Q(\phi,\nabla\phi,\nabla^{2}\phi)$
may be found in [4], Lemma C.2.
{}From this we see that $\Sig_{\phi}$ will have constant mean curvature
$\h_{\phi}\equiv1$ if and only if $N(\phi)=0$.
Moreover the linearized (or Jacobi) operator applied to $\phi$ is simply
$$
L\phi=\left.\frac{\del}{\del t}N(1+t\phi)\right|_{t=0}
=\Lap{\phi}+|\A|^{2}\phi. \tag 1.2
$$
Analyzing this operator, and a modified one introduced below,
allows us to establish the basic results concerning
the local structure of $\Mk$; in particular, the nonexistence of a
$L^{2}-$nullspace for $L$ at $\Sig\in\Mk$ guarantees that $\Mk$ is
locally a real analytic orbifold near $\Sig$.

We now state our main result concerning the local structure of $\Mk$.

\proclaim{Theorem 1.3} $\Mk$ is locally a finite dimensional real analytic
variety.  Moreover, if $\Sig\in\Mk$ satisfies the hypothesis:
$$
{\text{ If }}\phi\in L^{2}(\Sig) {\text { and }}
L\phi=0,{\text{ then }}\phi=0 \tag 1.4
$$
then in a neighborhood of $\Sig$, $\Mk$ is the quotient of a
real analytic manifold of dimension $3k-6$ by the finite isotropy subgroup
of $\Sig$ in the group of Euclidean motions,
for $k\geq 3$.
\endproclaim

The proof of Theorem 1.3 is presented in Theorems 3.1 and 4.1 below.

\smallskip

\noindent{\bf Remark 1.5.} For many CMC surfaces this isotropy
subgroup is trivial; however, there are numerous examples of CMC surfaces
for which this subgroup is nontrivial (see [2] and [4]).

\smallskip

\noindent{\bf Remark 1.6.} As noted above, there do not exist
CMC surfaces with only one end.  When $k=2$, any such surface is Delaunay and
these form a $3k-5=1$ dimensional family.  The proof of Theorem 1.3
recovers this fact.  It is interesting to note that the construction
given by Kapouleas has $3k-6$ continuous parameters (see [4], Remark 4.6),
however it is unknown whether or not these give rise to continuous
families of surfaces ([4], Remark 2.4).  Theorem 1.3 gives compelling
evidence that, in the case that a surface $\Sig$ constructed as in [4]
satisfies hypothesis (1.4), the $3k-6$ continuous parameters of Kapouleas
do yield continuous families of CMC surfaces.  In this case
the construction given in [4] would produce an open set in $\Mk$.

\smallskip

\noindent{\bf Remark 1.7.} Theorem 1.3 actually is valid in a more
general context.  The surface $\Sig$ need not be embedded but only
`weakly' or `Alexandrov' embedded, i\.e\. $\Sig=F(\partial\Omega)$ where
$F$ is an immersion into $\RR^{3}$ of a connected, open domain
$\Omega\subset\RR^{3}$.  For such surfaces
the Alexandrov reflection argument of [7] is still valid and Theorem 2.2
(see below) still holds.  This is the main geometric condition that a CMC
surface $\Sig$ needs to satisfy in order for us to establish Theorem 1.3.

\smallskip

This work arose in part from a previous study by Mazzeo, Pollack and
Uhlenbeck [8] of the solutions to the singular Yamabe problem on
$\SS^{n}\setminus\Lam$ where $\Lam=\{p_{1},\ldots,p_{k}\}$ is a set of
$k$ points.
The problem there was to understand the space $\ML$ of all metrics conformal
to the standard metric which are complete on $\SS^{n}\setminus\Lam$ and have
constant scalar curvature.  The existence of such metrics was established by
Schoen [11] and the main result of [8] is that $\ML$ is locally a real
analytic variety and is a real analytic manifold of dimension $k$ provided
a condition analogous to (1.4) is satisfied.
The proof of the ``good case'' of Theorem 1.3, i.e. when $L$ has no
$L^{2}-$nullspace is completely analogous to the proof of [8, Corollary 5.5],
and is sketched below in the proof of Theorem 3.1.  The proof in
[8] concerning the structure of $\ML$ as a real analytic variety does not
generalize to the constant mean curvature setting.  This forced us to find
a different proof, which is given in \S 4 below. This analysis also
applies to the singular Yamabe problem and provides a new, straightforward
proof that $\ML$ is always locally a real analytic variety.

Many authors have remarked on the strong analogy between constant
mean curvature surfaces and constant scalar curvature metrics.
In fact, most theorems in one subject have a counterpart in the other.  This
note reinforces this analogy, however a deep explanation for the relationship
between these two problems, one intrinsic and the other extrinsic, is still
missing.

\specialhead II. Background Material \endspecialhead

In this section we record some background material, principally results
from [7] and [8], which will be used in the proof of Theorem 1.3.

Given a unit vector $\bold{a}\in \RR^{3}$ and a smooth function
$\rho=\rho(t,\th)$, for $t\in\RR,\text{ and } \th\in [0,2\pi)$,
the {\it cylindrical graph} of $\rho$, about the axis
$\{t\bold{a}:t\in\RR\}\subset\RR^{3}$, is the image of the mapping
$$
\bold{F}(t,\theta)=t\bold{a}+\rho(t,\th)\o(\th),
$$
with $\o(\th)=\bold{b}\cos\th +\bold{c}\sin\th$, and $(\bold{a,b,c})$ a
positively oriented orthonormal frame.  The {\it steepness} of the graph
is measured by
$$
v=\sqrt{1+\rho_{t}^{2}+(\frac{\rho_{\th}}{\rho})^{2}}=\frac{1}{\nu\cdot\o},
$$
where $\nu$ is the unit outward normal.
As mentioned in the introduction, the Delaunay
surfaces refer to the (embedded) surfaces of revolution with constant mean
curvature equal to 1.  Such a surface, when expressed as a cylindrical graph,
is prescribed by its axis and a function $\rho=\rho(t)$,
which satisfies the differential equation
$$
\frac{\rho_{tt}}{v^{3}} - \frac{1}{\rho v} + 1 = 0. \tag 2.1
$$
The positive functions $\rho(t)$ satisfying (2.1) are all periodic
and may be distinguished by their minimum value $\e\in (0,1]$, which we
refer to as the Delaunay parameter of the corresponding surface $D_{\e}$.
This parameter is the minimum `neck size' of the surface.
The maximum bulge $\mu$ is the maximum of $\rho_{\e}$,
$1\leq\mu< 2$, where $\mu\rightarrow 2$ as $\e\rightarrow 0$.
The main result of [7] is that any complete, embedded CMC surface with
finite topology is asymptotically Delaunay.

\proclaim{Theorem 2.2 [7]} If $\Sig$ is a complete, embedded CMC surface of
finite topology,
then for each end $E_{j}$ of $\Sig$ there exists an embedded Delaunay surface
$D_{\e_{j}}$ in $\RR^3$ to which $E_{j}$ converges exponentially.  That is,
for $x\in\RR^{3}$ with $|x|$ sufficiently large, we may write $E_{j}$ and
$D_{\e_{j}}$ as cylindrical graphs over a fixed axis,
$\{t\bold{a}:t\in\RR\}\subset\RR^{3}$, of functions
$\rho_{E_{j}}$ and $\rho_{D_{\e_{j}}}$ respectively, with
$$
|\rho_{E_{j}} - \rho_{D_{\e_{j}}}|\leq C\, e^{-\lam_{j}t},
$$
as $t\rightarrow +\infty$, where $C$ and $\lam_{j}$ are positive constants,
and similar estimates hold for all derivatives of
$\rho_{E_{j}}$ and $\rho_{D_{\e_{j}}}$.
\endproclaim

An important ingredient in our analysis is an understanding of the Jacobi
fields on the Delaunay surfaces (i.e. solutions of $L\phi=0$ on $D_{\e}$).
By separating variables in the equation $L\phi=0$, one may
realize all the Jacobi fields in terms of the solutions of the corresponding
family of second order ODE's indexed by the eigenvalues $k^{2},\
k=0,1,2,\ldots$, of $-\del^2_{\th}$ on $\SS^{1}$.
There are six geometrically natural Jacobi fields on $D_{\e}$.  The first
two correspond to the infinitesimal translations along the Delaunay axis and
the infinitesimal  changes in the Delaunay parameter, and are denoted by
$\phi_{0,1}^{(\e)}$ and $\phi_{0,2}^{(\e)}$ respectively.  As functions of the
variables $(t,\th)$, these are the only Jacobi fields which are independent
of $\th$, i.e. they are the $k=0$ solutions.
The four solutions corresponding to the first
eigenvalue $k=1$ are the infinitesimal translations orthogonal to the Delaunay
axis and the infinitesimal rotations of the Delaunay axis about an
orthogonal axis.  These
are denoted by $\phi_{1,1,i}^{(\e)}$ and $\phi_{1,2,i}^{(\e)}$ respectively,
where $i=1,2$ denotes the two orthogonal axes.  Note that these functions
depend on the background Delaunay surface, and in particular on the Delaunay
parameter $\e$.
In both the $k=0$ and $k=1$ cases, the solutions corresponding to translations
are bounded and periodic, while the solutions corresponding to changes in the
Delaunay parameter and rotations of the axis are of linear growth.
As observed in [7], it is easy to see that all solutions corresponding to
the higher eigenvalues, $k\geq 2$, grow exponentially either as
$t\rightarrow\infty$ or $t\rightarrow-\infty$.

Using the fact that the CMC surface $\Sig$ has an asymptotic axis
corresponding to each end, we describe a neighborhood of
$\Sig$, containing an open set in $\Cal{M}_k$, as follows.
Some elements of $\Cal{M}_k$ near to
$\Sig$ have ends with axes in the same direction as those of $\Sig$. These
surfaces can be obtained as normal graphs over $\Sig$, as described
earlier. For other nearby surfaces, the directions of their axes may
differ by a small rotation from those of $\Sig$. Accordingly, we
form a $2k$ dimensional family $\Sig(\tau)$ of variations of $\Sig$.
Here $\tau = (\tau_1, \dots, \tau_k)$, where each $\tau_j = (\tau_{j,1},
\tau_{j,2})$ is a pair of real numbers describing a rotation of the
axis for the end $E_j$, pivoting near $\del E_j$, in the two
directions  orthogonal to this axis, which we index as $1,2$.
Recall that $\del E_j \subset
\del B_R(0)$ for $R$ sufficiently large. The exact transitions between
$\Sig \cap B_R(0)$ and the ends of each surface in this family are not
too important, provided it is done smoothly and `economically.'
However, note that each such surface has constant mean curvature except
in a small neighborhood of each $\del E_j$, and also that $\Sig(0) = \Sig$.
Since every CMC surface has asymptotic axes, by Theorem 2.2, we can
express any such surface near $\Sig$ as a normal graph over some
$\Sig(\tau)$ with $\tau$ small.

\head Weighted Sobolev Spaces and the Deficiency Subspace \endhead

The function spaces we use here are exponentially
weighted Sobolev spaces based on $L^2(\Sig)$; these are written
$H^s_{\ga}(\Sig)$, or just $H^s_{\ga}$, for
$\ga, s \in \RR$, $s>3$. The last condition ensures
that the spaces behave well under nonlinear operations, such as
taking products (e\.g\. terms like $\phi\nabla^{2}\phi$ which occur in
$Q(\phi,\nabla\phi,\nabla^{2}\phi)$).
To define $H^s_{\ga}(\Sig)$, decompose $\Sig$ into the
union of the ends $E_1, \dots, E_k$ and a compact piece $K$. Over $K$ an
element $h \in H^s_{\ga}$ restricts to an ordinary $H^s$ function.
Over $E_j$, $h = e^{\ga t}\tilde h$, where $\tilde h \in H^s([0,\infty)
\times \SS^{1}, dt\,d\th)$. Note that by Theorem 2.2 the measure here is
uniformly equivalent to the one induced by $\Sig$ for any fixed $\Sig \in
\Mk$.  When $\ga>0$, functions in $H^s_{\ga}(\Sig)$ are allowed to have
exponential growth, such as $e^{\ga^{\prime}t}$, for any $\ga^{\prime}<\ga$,
while functions in $H^s_{-\ga}(\Sig)$ must decay at least as fast as
$e^{-\ga t}$.  In particular, $H^s_{-\ga}(\Sig)\subset H^s(\Sig)\subset
H^s_{\ga}(\Sig)$.

In analogy with [8], given $\Sig$, a complete, embedded CMC surface with
$k$ ends, we use the Jacobi fields above to define a $6k-$dimensional linear
space, $W$, which we shall call the {\it{deficiency subspace}}.
Theorem 2.2 implies
there are $k$ axes and $k$ Delaunay surfaces $\{D_{\e_{j}}\}, j=1,\ldots,k$,
about those axes to which $\Sig$ converges as $|x|\rightarrow\infty$.
We define $W$ to be the linear span of the functions $\phi_{0,1}^{(\e_{j})}$,
$\phi_{0,2}^{(\e_{j})}$, $\phi_{1,1,i}^{(\e_{j})}$, and
$\phi_{1,2,i}^{(\e_{j})}$ for $i=1,2$ and $j=1,\ldots, k$, cut off to have
support outside of a ball $B_{R}(0)\in\RR^{3}$ for some $R$ sufficiently
large as in Theorem 2.2.  We note that $W\subset H^s_{\ga}(\Sig)$ for all
$\ga>0$ however $W\not\subset H^s(\Sig)$, since the Jacobi fields defining
$W$ are not in $L^2(\Sig)$.  Theorem 2.2 implies that the linearized operator
(1.2) on $\Sig$ is an asymptotically periodic operator on each end $E_j$.
In particular, on $E_j$ we have
$$
L = L_{\e_j} + e^{-\a t}F,   \tag 2.3
$$
where $F$ is a second order operator with coefficients bounded
in $\cin$ as functions of $(t,\th)$, $L_{\e_j}$ is the Jacobi operator
on $D_{\e_{j}}$  and $\a$ is a positive constant.

\head Fredholm Theory \endhead

Operators of the form (2.3) have arisen previously in problems in
geometry. In particular, we refer to work of C. Taubes [13] where some aspects
of the Fredholm theory in a similar asymptotically periodic setting were
developed.  In [8] a detailed analysis of certain mapping
properties for such operators was obtained.
The basic Fredholm result established there applies to
the linearization $L$ to give the following.
\proclaim{Proposition 2.4 [8]} There exists an infinite,
discrete set of numbers $\Gamma \subset \RR$ such that the bounded operator
$$
L:  H^{s+2}_{\gamma}(\Sig) \longrightarrow H^{s}_{\gamma}(\Sig) \tag 2.5
$$
is Fredholm for all values of the weight parameter $\gamma \notin\Gamma$.
In particular, $0 \in \Gamma$, so the map (2.5) is not Fredholm
on the ordinary unweighted Sobolev spaces, but
is Fredholm for all values of $\gamma$ sufficiently near, but not
equal to zero.
\endproclaim

The numbers $\ga_{j}\in \Gamma$ are the real parts of the `indicial roots'
of the model operators $L_{\e_j}$ on each end $E_{j}$.  The significance
of these indicial roots is explained, to some extent, in the discussion
of the relative index theorem below.  A more complete description of
$\Gamma$ and the proof of Proposition 2.4 is given in [8].

Understanding the local structure of the space $\Mk$ near $\Sig$ is closely
related to  determining when $L$ is actually injective or surjective,
and when not injective, understanding the elements in its nullspace.
In general this is a very difficult question; however, using duality and
the fact that $L$ is self-adjoint on $L^2$, i\.e\. when $\ga = 0$, we have
the following corollary.
\proclaim{Corollary 2.6} Suppose that $L$ has no $L^2-$nullspace.
Then for all $\d > 0$ sufficiently small
$$
\align
&L: H^{s+2}_{\delta} \longrightarrow H^s_{\delta} \quad \text{is
surjective},\\
&L: H^{s+2}_{-\delta} \longrightarrow H^s_{-\delta} \quad \text{is
injective}. \endalign
$$
\endproclaim

Corollary 2.6 states that if $L$ has no global $L^2-$nullspace, then we
can find a solution $w \in H^{s+2}_{\d}$ to the equation $Lw = f$ for every
$f \in H^s_{\d}$, whenever $\d > 0$. In particular, this holds whenever
$f \in H^s_{-\d}$.  Whenever $f$ decays at some exponential rate like this,
we expect the solution $w$ to be somewhat better behaved than
a general $H^{s+2}_{\d}$ function.  It is immediate
that $w$ is in $H^{s+2}_{\d}$ for any $\d>0$, however we can do even better.
This is the subject of what we will call the

\proclaim{Linear Decomposition Lemma 2.7 [8]}
Suppose $f \in H^s_{-\d}$ for some
$\d>0$ sufficiently small, and $w \in H^{s+2}_{\d}$ solves
$Lw = f$. Then $w \in H^{s+2}_{-\d} \oplus W$, i\.e\. $w$ may be decomposed
into a sum $v + \phi$, where $v\in H^{s+2}_{-\d}$ decays at the same
rate as $f$ and $\phi$ is in the deficiency subspace $W$.
\endproclaim

The space $W$ may also be thought of as the `parameter space' for the
linear problem, since its elements are the potentially occurring parameters
(disregarding the exponentially decreasing components) for Jacobi fields,
i\.e\. those $v$ such that $Lv=0$.  These functions are the subject of the
next subsection.

\head The Bounded Nullspace   \endhead

It remains to understand the nullspace of the Jacobi operator $L$.
When acting on $H^{s+2}_{\d}$ (for $\del>0$) the nullspace of $L$
will be the direct sum of the $L^2-$nullspace and what we call the
`bounded nullspace' $\BL$, defined by
$$
\BL = \{v \in H^{s+2}_{\d}: Lv = 0, v \notin H^{s+2}_{-\d} \}. \tag 2.8
$$
It follows immediately from the linear decomposition lemma (2.7) that
$\BL\subset H^{s+2}_{-\d} \oplus W$.
In addition one can show that (2.3) implies that the nullspace of
$L$ in  $L^2$ is the same as the nullspace of $L$ in $H^{s+2}_{-\d}$, i\.e\.
any solution $v\in L^2$ to $Lv=0$ decays exponentially at some uniform rate
(see Proposition 4.14, [8]).
Analyzing the $L^2-$nullspace of $L$ is a difficult problem which we do not
address here.  However, the dimension of $\BL$ may be calculated by means a
relative index theorem. We define the relative index for any
$\ga_{1}, \ga_{2}  \notin \Gamma$ by
$$
\relind(\ga_1,\ga_2) = \ind(\ga_1) - \ind(\ga_2).
$$
where
$$
\ind (\ga) = \dim\text{\,ker\,}\left. L\right|_{H^s_{\ga}}
- \dim\text{\,coker\,}\left. L \right|_{H^s_{\ga}}.
$$
{}From the fact that the adjoint of $L$ on
$H^s_{\ga}$ is $L$ on $H^{-s}_{-\ga}$, and using duality and
elliptic regularity, it is easy to see that for $\d>0$ sufficiently
small we have
$$
\relind(\d,-\d) =  2\dim(\BL).
$$
This reduces the problem of calculating $\dim(\BL)$ to that
of calculating $\relind(\d,-\d)$.  The latter may be calculated by using a
Relative Index Theorem of R. Melrose for metrics with asymptotically
cylindrical ends (see [10]).  The deformation $\Sig_{t}$ of $\Sig=\Sig_{0}$,
to a surface $\Sig_{1}$ with exactly cylindrical ends may be done in such a
way that the corresponding family of operators $L_{t}=\Lap_{\Sig_{t}} +
|{\A}_{\Sig_{t}}|^{2}$ are Fredholm on $H^s_{\d}$ and $H^s_{-\d}$
for every $0 \le t \le 1$.  This is a consequence of the fact that
the Jacobi fields corresponding to the $k>1$ eigenvalues grow or decay
exponentially on $D_{\e}$ for all $\e\in(0,1]$.  This means that
$\d\not\in\Gamma_{t}$ for any $0\leq t\leq 1$, where $\Gamma_{t}$
is the set of weights where $L_{t}$ fails to be Fredholm as in Proposition 2.4.
Under such deformations the relative index of the operators $L_{t}$ is
independent of $t$.

The relative index is computed in terms of asymptotic data on this deformed
surface. Since all ends of $\Sigma_1$ are cylindrical, the model operator
at each end is simply the linearization $L$ on the cylinder,
$\del_t^2 + \del_{\th}^2 + 1$. The asymptotic data we need to compute
is the sum of the  `ranks' of the poles along the real line of the
`indicial family' of this model operator, as defined in [10].
The indicial family is obtained by taking the Fourier transform in the
$t$ variable. Letting $\zeta$ be the dual variable, we obtain the
indicial family $ I(\zeta)\equiv\del_{\th}^2 + 1 - \zeta^2$.
This is a holomorphic family of elliptic operators whose
inverses, $I(\zeta)^{-1}$, form a meromorphic family of operators.
The poles of this family occur at the elements of
$$
\Cal P = \{\pm1, 0,\pm i\sqrt{k^2 - 1}\} {\text{ for }} k \ge 2.
$$
The only poles along the real line are at $\pm1$ and at $0$.
The poles at $\pm 1$ are simple, and have rank $1$.
This corresponds to the fact that these arise from the
eigenvalue  $\lam_{0}=0$ of $-\del^2_{\th}$, which is simple.
The pole at $0$ arises from
the multiplicity two eigenvalue $\lam_{1}=1$ of $-\del^2_{\th}$.
Since this pole
is double, the rank at $0$ counts each of these eigenvalues twice.
Hence, the rank at $0$ is $4$.
Thus the total rank, from all these poles, is $6$. Summing over all
$k$ ends, we conclude that the relative index is $6k$. Hence

\proclaim{Theorem 2.9} $\dim(\BL) = 3k$.
\endproclaim

The argument outlined above is exactly analogous to Theorem 4.24 of [8] and
we refer there for more details of the proof of Theorem 2.9.

\specialhead III.  The Smooth Case \endspecialhead

In this section we indicate how the results of \S 2 establish
that $\Mk$ is a  real analytic orbifold locally near $\Sig$
under the assumption that the linear operator $L$ has no $L^{2}-$nullspace.

\proclaim{Theorem 3.1} Suppose that $\Sig\in\Mk$ satisfies hypothesis (1.4).
Then there is an open set $\Cal U \subset \Mk$ containing $\Sig$, such that
$\Cal U$ is  the quotient of a
real analytic manifold of dimension $3k-6$ by the finite isotropy subgroup
of $\Sig$ in the group of Euclidean motions,
for $k\geq 3$.  When $k=2$, $\Cal U$ is  a
real analytic manifold of dimension $3k-5=1$; indeed ${\Cal M}_{2}\cong
(0,1]$.
\endproclaim

\demo{Proof}
We need to modify the nonlinear operator defined in (1.1) slightly
in order to accommodate normal variations off of the surfaces $\Sig(\tau)$,
with $\tau$ small, whose ends are rotations of those of $\Sig$. Consider
the operator
$$
N(\phi,\tau)\equiv 2(\h_{\phi,\tau}-1) = \Lap_{\Sig(\tau)}{\phi}+
|\A(\Sig(\tau))|^{2}\phi+
Q_{\Sig(\tau)}(\phi,\nabla\phi,\nabla^{2}\phi),
$$
where $\h_{\phi,\tau}$ is the mean curvature of the surface given by
the normal graph of $\phi$ over $\Sig(\tau)$.
We shall continue to refer to as the operator $N$.
The set $\Cal{U}$ is contained in the set of all solutions
of $N(\phi,\tau)=0$.

The main problem is to find spaces on which $N$ acts as a real
analytic function and on which its linearization $L$ is surjective.
Hypothesis (1.4) coupled with the linear decomposition lemma
implies that
$$
L: H^{s+2}_{-\d}\oplus W\rightarrow H^{s}_{-\d} \tag 3.2
$$
is a surjective map.  However, the nonlinear operator when applied to
elements of the form $(v+\psi)$ with $v\in H^{s+2}_{-\d}$ and $\psi\in W$
does not map into  $H^{s}_{-\d}$.  This may be remedied by recalling that
the elements of our parameter space $W$ are composed of functions
which are the derivatives of actual one parameter families of solutions to
$N_{\e_{j}}(v_{\eta}(t,\th))=0$ on each asymptotic Delaunay surface
$D_{\e_{j}}$ of $\Sig$. Given a Delaunay surface $D_{\e}$ there exist
functions $u_{\e+\eta}(t)$, $\eta\in (-\e,1-\e]$, which parameterize the
other Delaunay surfaces as deformations off of $D_{\e}$, i\.e\.
$$
D_{\e+\eta} = \{x+u_{\e+\eta}(t)\nu(x): x\in D_{\e}\}.
$$
In particular, on $D_{\e}$, $N(u_{\e+\eta}(t))=0$ and
$$
\phi_{0,2}^{(\e)} \equiv \left.\frac{d}{d\eta}u_{\e+\eta}(t)\right|_{\eta=0}
$$
satisfies $L\phi_{0,2}^{(\e)}=0$.  The curve $v_{\eta}=u_{\e+\eta}(t)$ is
a one parameter family of solutions to the nonlinear equation whose
tangent vector at $\eta=0$ is represented by the Jacobi field corresponding
to the infinitesimal changes in the Delaunay parameter on $D_{\e}$.
There are similarly defined curves corresponding to the other
Jacobi fields $\phi_{0,1}^{(\e)}$, and $\phi_{1,1,i}^{(\e)}$, for $i=1,2$.
The remaining Jacobi fields $\phi_{1,2,i}^{(\e)}$ in $W$ are obtained
as the derivatives of the family of surfaces $\Sig(\tau)$.
This allows us to define a corresponding 6k-dimensional manifold,
${\Bbb{W}}={\Bbb{W}}_{0}\times{\Cal{V}}$, where ${\Cal{V}}$ is a small
ball in the $\tau$ parameter space and $T_{0}{\Bbb{W}}_{0}$ is
the span of $\{\phi_{0,1}^{(\e)},\phi_{0,2}^{(\e)},\phi_{1,1,1}^{(\e)},
\phi_{1,1,2}^{(\e)}\}$.  ${\Bbb{W}}$ is composed of curves of
solutions which are annihilated by $N_{\e_{j}}$ on $E_{j}$.
By construction $T_{0}{\Bbb{W}}=W$.  Furthermore the nonlinear operator
$$
{N}: H^{s+2}_{-\d}\oplus {\Bbb{W}}\rightarrow H^{s}_{-\d}
$$
defined by $N(\eta_{1},\eta_{2},\tau)=N(\eta_{1}+\eta_{2},\tau)$
is a real analytic map with $L$ as its surjective linearization.

To verify the last claim we proceed as follows. The main step is
to write $N(\eta_{1}+\eta_{2},\tau)$ in such a way that it obviously
belongs to $H^{s}_{-\d}$.   On each end $E_{j}$ we may
write $N$ as an exponentially small perturbation of the corresponding operator
on the asymptotic Delaunay surface, i\.e\. $N=N_{\e}+e^{-\a t}\hat{Q}$.
We then use a common formula for the remainder in Taylor's theorem
to obtain
$$
N(\eta_{1}+\eta_{2},\tau) = N_{\e}(\eta_{2},\tau) +
e^{-\a t}{\hat{Q}}(\eta_{1}+\eta_{2},\tau) +
\left[\int_0^1 N^{\prime}(s\eta_{1}+\eta_{2},\tau)\,ds\right]\eta_{1}.
$$
Since $N_{\e}(\eta_{2},\tau) = 0$ and
$N^{\prime}(s\eta_{1}+\eta_{2},\tau)$ is the linearization $\Lap+|\A|^{2}$
for the surface given as the normal
graph of $s\eta_{1}+\eta_{2}$ over $\Sig(\tau)$,
every term on the right here is in $H^{s}_{-\d}$
so long as $\d<\a$. Furthermore, every term is real analytic
in $(\eta_{1}, \eta_{2},\tau)$.
Real analyticity of $N$ for functions supported
in the interior, away from the ends, follows directly from the
real analyticity of the mean curvature equation.

We now apply the real analytic implicit function theorem to conclude
that there is a set $\Cal U_{0}\subset  H^{s+2}_{-\d}\oplus {\Bbb{W}}$,
consisting of those functions annihilated by $N$ on  $\Sig$, such that
$\Cal U_{0}$ is a finite dimensional real analytic manifold.  The dimension
of $\Cal U_{0}$ is simply the dimension of the nullspace $L$ in (3.2).  This
is precisely the bounded nullspace $\BL$, whose dimension is $3k$ as shown in
Theorem 2.9.  The open set $\Cal U \subset \Mk$ is obtained by identifying
those functions in $\Cal U_{0}$ which arise from the (real analytic) space of
isometries of $\RR^3$.  The dimension
of $\Cal U$ is the dimension of the space of solutions to the partial
differential equation $N(\phi,\tau)=0$ minus the dimension of
${\text{Isom}}(\RR^{3})$ plus the dimension of ${\text{Iso}}(\Sig)$,
the isotropy subgroup  of $\Sig$.   When $k\geq 3$, ${\text{Iso}}(\Sig)$
is finite (and possibly trivial), and when $k=2$ it consists of a subgroup of
rotations (isomorphic to $\SS^{1}$) and a subgroup (discrete unless
$\e=1$ and $\Sig$ is the cylinder) of
translations of $\Sig=D_{\e}$ along the Delaunay axis. Thus
the dimension of $\Cal U$ is $ 3k-6$ when $k\geq 3$,
and $3k-5\equiv 1$ when $k=2$.
When $k\geq 3$ the existence of a finite isotropy subgroup of $\Sig$
implies that the action of ${\text{Isom}}(\RR^{3})$ on $\Cal U_{0}$
is not free and thus $\Cal U$ is a real analytic orbifold.
This completes the proof of Theorem 3.1
\enddemo

\smallskip

\noindent{\bf{Remark 3.3.}} Every element $\hat\Sig\in \Cal{U}\subset\Mk$,
may be represented by some $(\eta,\tau)=(\eta_{1},\eta_{2},\tau)\in
H^{s+2}_{-\d}\oplus {\Bbb{W}}$.  The function $\eta_{2}\in {\Bbb{W}}_{0}$
consists of some combination of changes in the asymptotic Delaunay parameters,
and asymptotic translations along and orthogonal to the axes of each end.
$\tau$ corresponds to rotations of these axes.
Determining precisely which
of these $6k$ parameters occur in the deformations of a given surface
$\Sig$ is an interesting geometric problem.  The
deformations which can occur as global solutions of the equation
$N(\phi,\tau)=0$ correspond, on the linear level, to the bounded nullspace
$\BL$.  We cannot determine $\BL$ explicitly, however we remark here
that if $W\cong\RR^{6k}$ is endowed with a natural symplectic structure,
then $\BL\subset W$ is a Lagrangian subspace.  We refer to ([8], Lemma 7.5)
for the proof of this fact.

\specialhead IV. The General Case \endspecialhead

In this section we establish the main result about the local structure
of $\Mk$ in the general  case, when the linearization $L$ may have
an $L^2-$nullspace.

\proclaim{Theorem 4.1} The space $\Mk$ is locally a finite dimensional
real analytic variety.
\endproclaim

\demo{Proof}  Given $\Sig\in\Mk$, we here consider the possibility
that the Jacobi operator $L$ is not injective when acting on the
space $H^{s+2}_{-\d}$ for any $0<\d<\d_{0}$.  (Note that the alternative
``good case'' has been covered in \S3.)
Proposition (2.4) implies that the bounded operator
$$
L: H^{s+2}_{+\d}\rightarrow H^{s}_{+\d}
$$
is Fredholm and, by duality, has a finite dimensional cokernel which we
identify with $K_{(1)}\equiv K_{-\d}$, the kernel of $L$ restricted to
$H^{s+2}_{-\d}$.  Consider the surjective Fredholm
operator
$$
\Cal{L}: H^{s+2}_{+\d}\oplus K_{(1)}\rightarrow H^{s}_{+\d}
$$
defined by
$\Cal{L}(\eta,\phi)=L(\eta)+\phi \in H^{s}_{+\d}$.

Since $H^{s}_{-\d}\subset H^{s}_{\d}$, we may now apply the
Linear Decomposition Lemma (2.7) to conclude that the operator
$$
\Cal{L}: H^{s+2}_{-\d}\oplus W\oplus K_{(1)}\rightarrow H^{s}_{-\d} \tag 4.2
$$
defined by $\Cal{L}(\eta_{1},\eta_{2},\phi)=
L(\eta_{1}+\eta_{2})+\phi \in H^{s}_{-\d}$ is surjective.
Recall that elements of our parameter space $W$ are composed of functions
which are the Jacobi fields corresponding to one parameters families of
Delaunay surfaces (corresponding to the deformations of
the asymptotic Delaunay surfaces of $\Sig$.)
This allowed us to define a corresponding $6k$ dimensional manifold,
${\Bbb{W}}$, annihilated by the model operators $N_{\e_{j}}$ on $D_{\e_{j}}$
such that $T_{0}{\Bbb{W}}=W$.  Thus the nonlinear operator
$$
\Cal{N}: H^{s+2}_{-\d}\oplus {\Bbb{W}}\oplus K_{(1)}\rightarrow H^{s}_{-\d}
$$
defined by $\Cal{N}(\eta_{1},\eta_{2},\tau,\phi)=
N(\eta_{1}+\eta_{2},\tau)+\phi$
is a real analytic map with $\Cal{L}$ as its surjective linearization.
We are interested in the space
$$
\aligned
\Cal{S}&=\{\tilde{\phi}\in H^{s+2}_{-\d}\oplus {\Bbb{W}}:
N({\tilde{\phi}})=0\}\cr
       &=\{(\tilde{\phi},\phi_{0})\in H^{s+2}_{-\d}\oplus {\Bbb{W}}\oplus
K_{(1)} : \Cal{N}(\tilde{\phi},\phi_{0})=\phi_{0}\}.
\endaligned
$$
Since $K_{(1)}\subset H^{s}_{-\d}$ is a linear subspace (in particular a
real analytic submanifold), the implicit function theorem shows that
$$
\Cal{T}=\{(\tilde{\phi},\phi_{0})\in H^{s+2}_{-\d}\oplus {\Bbb{W}}\oplus
K_{(1)}:\Cal{N}(\tilde{\phi},\phi_{0})\in K_{(1)}\}
$$
is a finite dimensional real analytic manifold.
In fact, $\Cal{T}=\ker(\Pi^{\perp}\circ\Cal{N})$, where
$\Pi:H^{s}_{-\d}\rightarrow K_{(1)}$ is the orthogonal projection and
$\Pi^{\perp}= {\text{I}}-\Pi$. The linearization of $\Pi^{\perp}\circ\Cal{N}$
is $\Pi^{\perp}\circ\Cal{L}$, and by (4.2) this is surjective.
The nullspace of $\Pi^{\perp}\circ\Cal{L}$ is the set of
$(\eta_{1},\eta_{2},\phi)$ such that $\Cal{L}(\eta_{1},\eta_{2},\phi)\subset
K_{(1)}$, since $\Pi^{\perp}$ annihilates any $\phi\in K_{(1)}$.
If $\Cal{L}(\eta_{1},\eta_{2},\phi)\in K_{(1)}$ then $\phi\in K_{(1)}$
is arbitrary and $L(\eta_{1}+\eta_{2})=0$, since the range of $L$ is
orthogonal to $K_{(1)}$ in $H^{s}_{-\d}$.
This implies that $\eta_{1}+\eta_{2}\in
K_{(2)}\oplus B$, where $K_{(2)}\equiv K_{-\d}$, is a second copy of the
kernel of $L$ acting on $H^{s+2}_{-\d}$, and $B$ is the bounded nullspace.
$\Cal{T}$ may thus be represented locally as a graph in
$\TT=H^{s+2}_{-\d}\oplus {\Bbb{W}}\oplus K_{(1)}$,
i.e. there exists a real analytic map
$$
\Psi:K_{(2)}\oplus B\oplus K_{(1)}\rightarrow
(H_{-\d}^{s+2}\ominus K_{(2)})\oplus (W\ominus B)
$$
such that
$$
\Cal{T}=\{\phi+\Psi(\phi)\in \TT : \phi\in K_{(2)}\oplus B\oplus K_{(1)}\}.
$$
Clearly $\Cal{S}\subset\Cal{T}$, in fact writing
$\phi=(\phi_{1},\phi_{2},\phi_{0})$, we have
$$
\Cal{S}=\{\phi+\Psi(\phi)\in \TT :\Cal{N}(\phi)=\phi_{0},
\phi\in K_{(2)}\oplus B\oplus K_{(1)}\}.
$$
Equivalently,
$$
\aligned
\Cal{S}=&\{\phi+\Psi(\phi)\in \TT :N(\phi_{1}+\phi_{2}+\Psi(\phi))=0,
\phi\in K_{(2)}\oplus B\oplus K_{(1)}\}\cr
=&\{\phi\in K_{(2)}\oplus B\oplus K_{(1)}:N\circ(\tilde{I}+\Psi)(\phi)=0\},
\endaligned
$$
where $\tilde{I}(\phi)=(\phi_{1},\phi_{2})$, is the projection onto the
first two factors.
Finally, since for any $\phi\in K_{(2)}\oplus B\oplus K_{(1)}$,
$\Cal{N}\circ(I+\Psi)(\phi)=N\circ(\tilde{I}+\Psi)(\phi)+\phi_{0}\in K_{(1)}$,
we see that $N\circ(\tilde{I}+\Psi)(\phi)\in K_{(1)}$.
Thus we have represented locally the space of solutions to our
nonlinear equation $N(v,\tau)=0$, with $v$ and $\tau$ small, as the zero set of
$$
N\circ(\tilde{I}+\Psi):K_{(2)}\oplus B\oplus K_{(1)}\rightarrow K_{(1)}
$$
which is a real analytic map between finite dimensional spaces.
As before, the space $\Mk$ is obtained from this space by identifying those
functions which arise from an isometry of $\RR^{3}$.
This completes the proof of Theorem 4.1.

\enddemo

\smallskip

\noindent{\bf{Remark 4.3.}}
The technique of the proof of Theorem 4.1 is adapted from
an argument of L. Simon ([12], Theorem 3), where a similar result is
established in the course of his extending the \L ojasiewicz inequality
to infinite dimensions.
This simple technique for establishing real analyticity for the solution
space of a real analytic elliptic operator
is well known in the setting of compact manifolds, and is often
referred to as the `Liapunov-Schmidt reduction'.

\specialhead V. CMC Hypersurfaces in Hyperbolic Space \endspecialhead

In this section we sketch how the results above apply to a class
of CMC hypersurfaces in $\HH^{n+1}$.  Generalizations of the Delaunay surfaces
to hypersurfaces in $\HH^{n+1}$ (and $\SS^{n+1}$) were given
by W.Y.~Hsiang [3]. He showed that there exists a
one parameter family of embedded, periodic CMC surfaces of revolution about
a geodesic line $\ga\in\HH^{n+1}$ which interpolate between a cylinder and
a string of geodesic spheres.

In [6], an analogue of Theorem 2.2 was established.
In dimension $n+1=3$, any complete surface in $\HH^{3}$ with constant mean
curvature greater than $2=n$ (that of a horosphere) has $k\geq2$ ends,
each of which is cylindrically bounded.  In particular any such surface
$\Sig$ is conformally a compact Riemann surface with finitely many punctures
(as is the case [7] for CMC surfaces in $\RR^3$), and the limit set of
$\Sig\subset\HH^{3}$ is a set of $k$ points in the sphere at infinity
$S_{\infty}$.  Any such surface is asymptotic to a generalized Delaunay
surface in the sense of Theorem 2.2.  When $n+1>3$ this is also true provided
one makes the {\it assumption} that the ends of $\Sig$ are
cylindrically bounded, i\.e\. remain within a bounded distance of some
geodesic ray.   For CMC hypersurfaces $\Sig\subset \HH^{n+1}$ cylindrical
boundedness of the ends implies that ${\bold{H}}>n$

We let $\Lambda = \{p_1, \cdots, p_k\}\subset S_{\infty}\cong \SS^{n}$
denote the limit set of a hypersurface $\Sig\subset\HH^{n+1}$ with constant
mean curvature ${\text{\bf H}}>n$. The space $\ML$ is defined to be
the set of all such hypersurfaces having $\Lam$ as their limit set,
where any hypersurfaces which differ by an isometry of $\HH^{n+1}$
are identified.

The local structure of $\ML$ is given by the following theorem.

\proclaim{Theorem 5.1} $\ML$ is locally a finite dimensional real analytic
variety.  Moreover, if $\Sig\in\ML$ satisfies hypothesis (1.4)
then in a neighborhood of $\Sig$, $\ML$ is the quotient of a
real analytic manifold of dimension $k$ for $k\geq 3$
(and of dimension 1 for $k=2$) by the finite isotropy subgroup
of $\Sig$ in the group of hyperbolic  motions, i\.e\.  $\ML$ is
a real analytic orbifold.
\endproclaim
The operator $L$ in Theorem 5.1 is now the linearization of the
quasilinear equation which specifies the mean curvature of a surface
$\Sig_{\phi}$, written as a graph over $\Sigma\in\HH^{n+1}$, to be
identically ${\text{\bf H}}$.
The results of [6] discussed above imply that this operator is an
asymptotically periodic operator as in (2.3).  In particular,
analogues of Proposition 2.4, Corollary 2.6, Lemma 2.7 and Theorem 2.9 hold.
The distinction is that the dimension of the bounded nullspace
is $k$ as opposed to $3k$.  This is due to the fact that the Jacobi fields
arising from the rotations and parabolic translations fixing either end of
the asymptotic Delaunay
axes about the $n$ orthogonal axes have exponential growth or decay and
hence do not contribute to the relative index.  The proof of Theorem 5.1
is identical to the proof given above for Theorems 3.1 and 41.

One may note that the dimension $k$ of the space $\ML$ is the same as
the dimension given in [8] for the solutions to the singular Yamabe
problem on $\SS^{n}\setminus\Lam$.  In the hyperbolic setting fixing the
limit set is precisely analogous to fixing the singular set
$\Lam\subset\SS^{n}$ in the singular Yamabe problem.  It was suggested to
us by Jonathan Poritz that one might try to recover this dimension count
for the moduli space of CMC surfaces in $\RR^{3}$ by specifying that the
asymptotic axes remain fixed as rays in $\RR^{3}$.  In other words,
we define a new space $\Mk^{0}\subset\Mk$, whose elements are those
surfaces in $\Mk$ which have the same asymptotic Delaunay axes as $\Sig$.
We conjecture that,
in the case when $\Sig$ satisfies hypothesis (1.4), $\Mk^{0}$ is also a real
analytic manifold of dimension $k$.  The point is that those Jacobi fields
arising from the rotations and translations of the asymptotic Delaunay
axes about the orthogonal axes should be ruled out since they correspond to
one parameter families of solutions about the asymptotic Delaunay surfaces
which change the Delaunay axis.  It is not clear how
this may be done in the proof of Theorem 1.3.

\vfill\eject

\enddocument

\bye